\documentclass{article}
\usepackage{amsmath}
\usepackage{amssymb}
\usepackage{amsthm}
\usepackage[shortlabels]{enumitem}
\usepackage{breqn}
\usepackage{authblk}

\newtheorem{theorem}{Theorem}

\newtheorem{corollary}[theorem]{Corollary}

\newtheorem{definition}[theorem]{Definition}

\newtheorem{lemma}[theorem]{Lemma}

\newtheorem*{theorem*}{Theorem}

\newcommand{\K}{\mathcal{K}}

\newcommand{\Len}{\mathcal{L}}
\newcommand{\Cons}{\mathcal{C}}
\newcommand{\A}{\mathbf{A}}
\newcommand{\B}{\mathbf{B}}

\newcommand{\N}{\mathbb{N}}
\newcommand{\Mod}{\mathrm{Mod}}

\title{Herbrand's Theorem: a short statement and a model-theoretic proof}
\author{Mariana Badano}
\affil{Facultad de Matemática, Astronomía, Física y Computación (FAMAF), Universidad Nacional de Córdoba, Argentina. \\ \texttt{marianabadano@unc.edu.ar}}
\date{}
\begin{document}
\maketitle
\begin{abstract}Herbrand's Theorem is a fundamental result in mathematical logic which provides a reduction of first-order formulas satisfied by a universal class to formulas free of existential quantifiers. In this work, a simpler and self-contained formulation of Herbrand's Theorem is presented, along with a model-theoretic proof of its general version.\end{abstract}

\noindent \textbf{Keywords:} Herbrand's theorem, Universal classes, Quasivarieties, Sound functions.

\noindent \textbf{Mathematics Subject Classification (2020):} Primary 03C07; Secondary 03C40, 03F03, 08C15.

\section{Introduction}
Herbrand's Theorem was first stated by Herbrand in his 1930 thesis \cite{herb}. It is a fundamental result in mathematical logic and has served as the theoretical foundation for many proof procedures.

The theorem commonly referred to as ``Herbrand's Theorem'' in modern logic courses is typically a weaker version, restricted to formulas which, in prenex form, contain only existential quantifiers. In contrast, the full version of Herbrand's Theorem applies to all first-order formulas. The simpler version has a well-known model-theoretic proof, which proceeds via a simple expansion of the language with constants, preservation under substructures, and an application of the Compactness Theorem, and is usually presented as an exercise in textbooks such as \cite{Chang} and \cite{Hodges}. The original proof of the general theorem is a highly complex syntactic proof. This approach has been extensively studied in the detailed survey by Buss \cite{Buss}. Furthermore, the absence of a direct and uniform statement for the general version makes the theorem virtually inaccessible to many researchers.
In this work, we provide a semantic proof of the theorem; more significantly, we present a simple statement of it.

To illustrate the inherent difficulties in formalizing the result, we begin by considering how the general version of the theorem could informally be stated:
Let $\Sigma$ be a set of universal sentences and $O$ an open formula. Let
\[\Upsilon:= \exists y_1\forall z_1\dots \exists y_k\forall z_k \; O(y_1,z_1,\dots,y_k,z_k).\]
Then, the following are equivalent:
\begin{enumerate}[label=\textup{(\arabic*)}]
\item $\Sigma\vdash \Upsilon$;
\item There are terms $t_{ij}$ and ``\textit{well-chosen}'' variables $x_{ij}$, with $ 1\leq i \leq N$ and $1\leq j\leq k$, such that $\Sigma\vdash \bigvee_{i=1}^N O(t_{i1},x_{i1},\dots,t_{ik},x_{ik})$, for some $N>0$.
\end{enumerate}

The problem with this attempt at a statement is that it is not at all clear what is meant by ``well-chosen'' variables. Quite an evident requirement is that any variable $x_{ij}$ should not occur in a $t_{im}$ with $m\leq j$; but this condition is not strong enough. For example, let $\Sigma:=\{\forall x \ (x{=}0\lor x{=}1)\}$ in the language $\{0,1\}$, where $0$ and $1$ are constant symbols. By taking $O(y_1,z_1):= z_1{=}y_1$, $t_{11}:=0$, $t_{21}:=1$ and $x_{11},x_{21}:=x$, we have $\Sigma\vdash \bigvee_{i=1}^2 O(t_{i1},x_{i1})$, but $\Sigma\nvdash \exists y_1 \forall z_1 \ O(y_1,z_1)$; that is, $(2)$ does not imply $(1)$. We could then propose that all $x_{ij}$ be distinct; however, this condition is too strong, as it is not true that $(1)$ always implies $(2)$ with that restriction. This highlights the difficulty of finding a set of conditions on the variables that preserves both directions of the equivalence, a challenge that has also affected previous attempts to formalize the result.

In 1956, {\L}o\'{s}, Mostowski, and Rasiowa in \cite{Most} presented a statement of this kind, proposed a notion of ``well-chosen'' variables, and gave an alternative to the original proof via algebraic methods. In 1961, after Feferman found an error in the proof of the direction $(1)\Rightarrow (2)$, the authors corrected it in \cite{Most-Fefe}. We pinpoint two crucial issues in those works. First, the statement of Herbrand's Theorem is not general. The theorem is stated for a formula with at most five nested quantifiers: three existential and two universal. The complexity in stating the theorem is determined by the number of universal quantifier levels. Although the authors claim that this case is general enough to illustrate the proof for any number of quantifiers, the general statement is never provided.
Secondly, and more importantly, the statement of Herbrand's Theorem presented is false. The notion of ``well-chosen'' proposed by the authors is too weak and does not guarantee $(2)\Rightarrow (1)$. In Section 2, we present a simple example where the implication $(2)\Rightarrow (1)$ of the theorem fails for the theory of linear orders.

In order to state the theorem properly, in Section 3, we define a precise type of functions, called \textit{sound functions}, that code tuples of terms to ensure the correct dependency of variables. We prove the existence of sound functions for countable languages by constructing one via a G\"{o}del numbering.
We give a semantic proof of $(1)\Rightarrow (2)$ by means of model-theoretic methods closely related to those used for the simple case of an existential formula. By Compactness, we assume, without loss of generality, that the language $\Len$ involved is finite and consider any sound function. We then expand the language with a countable set of new constants and consider the class of structures satisfying $\Sigma$ in the expanded language $\Len^*$. For an arbitrary $\Len^*$-structure $\A$ satisfying $\Sigma$, we consider the substructure $\B$ generated by the empty set. Since $\Sigma$ is a set of universal sentences, it is preserved under substructures; hence $\B\models\Sigma$ and therefore $\B\models\Upsilon$. This implies the existence of ground terms that witness the existential quantifiers, and we select constants indexed by a sound function applied to those terms as instances of the universally quantified variables. This gives a ground instance of $O$ satisfied by $\B$ and therefore satisfied by $\A$. This gives us a set of ground instances of $O$ such that every $\Len^*$-structure satisfying $\Sigma$ satisfies at least one of them. Then, by the Compactness Theorem, we obtain a finite disjunction of them. Finally, by replacing the new constants with variables, we obtain $(2)$.
For the proof of $(2)\Rightarrow (1)$, we proceed by syntactic methods; in fact, we prove that the finite disjunction described in (2), when indexed by any sound function, logically implies $\Upsilon$.

We believe that our model-theoretic approach yields a proof that is significantly simpler and more transparent than earlier approaches. Overall, the main contribution of this article is the applicability of the formulation of Herbrand's Theorem we present. We hope it will facilitate the use of the theorem within various areas where the definability of objects by first-order formulas with an arbitrary number of quantifiers is employed. This occurs, for example, in the theory of central elements (see \cite{bava2020} and \cite{Vaggione1999}), where the definability of factor congruences can be achieved using formulas which, although preserved by products and direct factors, may include arbitrary quantifiers. More generally, in the final section, we present two results obtained by the direct application of this presentation of the theorem, which show the advantage of a general statement, independent of the number of quantifiers.

\section{A note on a previous formulation of Herbrand's Theorem}
In this section, we will give a counterexample to the statement of Herbrand's Theorem proposed by {\L}o\'{s}, Mostowski and Rasiowa. First, we present a summary of the framework developed in \cite{Most}.

Consider a theory $T$ axiomatized by a set of universal axioms over the language $\Len$. (To ensure that the set of terms can be arranged in a countable sequence as the authors assume, the language $\Len$ must be taken to be countable. This does not restrict generality, as any specific deduction $T\vdash \Upsilon$ takes place within a finite sublanguage.)
Let $y_1, y_2,\dots$ be a countable sequence of variables and let $t_1,t_2,\dots$ be an enumeration of the terms in $\Len$.
Consider an open formula $O(y_1,y_2,y_3,y_4,y_5)$ of $\Len$ and let
\[\Upsilon:=\exists y_1\forall y_2\exists y_3\forall y_4\exists y_5 \ O(y_1,y_2,y_3,y_4,y_5).\]
Now, let $\mathbf{k}:\N\to \N$ and $\mathbf{l}:\N\times \N\to \N$ be functions satisfying the following properties:
\begin{enumerate}[label=(\roman*)]
\item $\mathbf{k}$ is a strictly increasing function;
\item $\mathbf{k}(n)>5$, for every $n\in \N$;
\item $y_{\mathbf{k}(n)}$ does not occur in $t_n$;
\item $\mathbf{l}(m,n)\neq \mathbf{k}(r)$, for every $m,n,r\in \N$;
\item if $\mathbf{l}(m,n)=\mathbf{l}(m',n')$ then $m=m'$ and $n=n'$;
\item $y_{\mathbf{l}(m,n)}$ occurs in neither $t_n$ nor $t_m$.
\end{enumerate}
For every $n\in \N$, let
\[H_n=\bigvee_{j=1}^n\bigvee_{l=1}^n\bigvee_{h=1}^n O(t_j,y_{\mathbf{k}(j)},t_l,y_{\mathbf{l}(j,l)}, t_h).\]
The statement of Herbrand's Theorem presented there is as follows.

\begin{theorem*}[{\L}o\'{s}, Mostowski and Rasiowa \cite{Most}]
The following conditions are equivalent:
\begin{enumerate}[label=\textup{(\arabic*)}]
\item $\Upsilon$ is provable in $T$;
\item there is an $n$ such that $H_n$ is provable in $T$.
\end{enumerate}
\end{theorem*}

This statement of Herbrand's Theorem is false, even for the case of five quantifiers. This issue is unrelated to Feferman's observation of an error in the proof, which appears in the implication $(1)\Rightarrow (2)$ of the theorem and was corrected in a later work \cite{Most-Fefe}. The part of the theorem we claim to be false is the implication $(2)\Rightarrow(1)$. The problem lies in the definition of the functions $\mathbf{k}$ and $\mathbf{l}$. Since $\mathbf{k}$ and $\mathbf{l}$ are required to be injective, they prevent counterexamples such as the one presented in our Introduction, but these functions do not prevent cycles. Although it is required that $y_{\mathbf{k}(n)}$ does not occur in $t_n$, situations of the form ``$y_{\mathbf{k}(n)}$ occurs in $t_m$ and $y_{\mathbf{k}(m)}$ occurs in $t_n$'' are admitted. This interweaving of variables is also permitted when both functions $\mathbf{k}$ and $\mathbf{l}$ are considered. We present a simple example where the stated theorem fails.

Let $\Sigma$ be a set of universal axioms for the theory of linear orders over the language $\Len = \{\leq\}$. Let $O(y_1,y_2):=y_1\leq y_2$ and $\Upsilon:=\exists y_1 \forall y_2 \ O(y_1, y_2)$. This case is effectively a specialization of the five-quantifier schema where the matrix is independent of the additional variables $y_3, y_4, y_5$, and the extra quantifiers are redundant. In this simplified setting, the formula $H_n$ reduces to
\[H_n=\bigvee_{j=1}^n O(t_j,y_{\mathbf{k}(j)}).\]
Consider an enumeration of the terms such that $t_1 = y_8$ and $t_2 = y_6$, and a function $\mathbf{k}$ such that $\mathbf{k}(1) = 6$ and $\mathbf{k}(2) = 8$. Note that $\mathbf{k}(1), \mathbf{k}(2) > 5$, $y_{\mathbf{k}(1)} = y_6$ does not occur in $t_1$, $y_{\mathbf{k}(2)} = y_8$ does not occur in $t_2$, and $\mathbf{k}(1)<\mathbf{k}(2)$. Therefore, the conditions imposed on the function $\mathbf{k}$ are satisfied locally, and we can extend $\mathbf{k}$ and define $\mathbf{l}$ on $\N$ to satisfy the global requirements.
Since $\Sigma$ axiomatizes linear orders, we have
\[ \Sigma \vdash y_6 \leq y_8 \lor y_8 \leq y_6, \]
which translates to:
\[ \Sigma \vdash t_2 \leq y_{\mathbf{k}(2)} \lor t_1 \leq y_{\mathbf{k}(1)}, \]
that is, $\Sigma \vdash H_2$.
By $(2)\Rightarrow (1)$ of the theorem, we infer
\[\Sigma\vdash\exists y_1 \forall y_2 \ y_1\leq y_2,\]
which is false, since not every linear order has a least element.

It is worth noting that, even if correctly defined, the functions $\mathbf{k}$ and $\mathbf{l}$ from the characterization by {\L}o\'{s}, Mostowski, and Rasiowa only describe the choice of variables for two levels of universal quantifiers. In the framework of our informal description, this corresponds to the case $k=2$, where $\mathbf{k}$ and $\mathbf{l}$ would match the restrictions to $T(X)$ and $T^2(X)$, respectively, of the general notion of a sound function introduced in the next section.

\section{Main Theorem}
For a set $\Sigma$ of sentences of a language $\Len$ and an $\Len$-structure $\A$, we use the notation $\A\models \Sigma$ to state that $\A\models\varphi$ for every $\varphi\in \Sigma$. We denote by $\Mod_{\Len}(\Sigma)$ the class of $\Len$-structures $\A$ such that $\A\models \Sigma$. We say that $\Sigma$ is a set of axioms for a class $\K$ if $\K=\Mod_{\Len}(\Sigma)$. We say that $\K$ is a universal class if it has a set of axioms consisting of universal sentences. For a sentence $\varphi$ and a class $\K$, the notation $\K\models \varphi$ states that, for every $\A\in \K$, $\A\models \varphi$. If $\Sigma$ is a set of sentences, the notation $\K\models \Sigma$ means $\K\models \varphi$ for every $\varphi\in \Sigma$. If $\varphi$ is a formula with free variables among $x_1,\dots,x_n$, the notation $\K \models \varphi$ means $\K\models \forall x_1\dots\forall x_n \; \varphi(x_1,\dots,x_n)$.

Let $X=\{x_i: i\in \N\}$, $Y=\{y_i: i\in \N\}$, and $Z=\{z_i: i\in \N\}$ be three pairwise disjoint sets of variables. For an arbitrary first-order language $\Len$, let $T_\Len(X)$ be the set of terms over $X$. We denote the set of ground terms of $\Len$ by $T_\Len(\emptyset)$. For a finite, non-empty subset $\{x_{1},\dots, x_{n}\}\subseteq X$, we denote the set of terms over $\{x_{1},\dots, x_{n}\}$ by $T_\Len(x_{1},\dots, x_{n})$. In each of these cases, we drop the subscript $\Len$ when there's no ambiguity about the language.

We consider first-order formulas in prenex normal form with an alternating quantifier structure that starts with an existential quantifier, specifically, formulas of the form
\[
\exists y_1\forall z_1\dots\exists y_k\forall z_k\;\varphi(y_1,z_1,\dots,y_k,z_k)
\]
with $k\geq 1$ and $\varphi$ an open formula.
The reader can easily verify that any first-order formula is logically equivalent to a formula of the form above (once the formula is in prenex form, we simply add the necessary extra dummy quantifiers).

To state our version of Herbrand's Theorem, we introduce a special type of function. This formalizes the notion of ``well-chosen'' variables mentioned in the Introduction.
\begin{definition}\label{sound} \textup{ For a given set of $k$-tuples of terms $ \mathcal{T}\subseteq T^{k}(X)$, we say that $\kappa: \{\langle t_{1},\dots,t_{j}\rangle: \vec{t}\in \mathcal{T}\text{ and }1\leq j\leq k\}\to \N$ is a \textit{sound function} \textit{for} $\mathcal{T}$ if }
\begin{enumerate}[\textup{(\roman*)}]
\item \label{itiii} $\kappa$\textup{ is injective,}
\item \label{iti} $\kappa( t_{1},\dots,t_{j})<\kappa( t_{1},\dots,t_{j+1})$\textup{ for $j<k$, and }
\item \label{itii}\textup{if $x_{m}$ occurs in $t_{j}$, then $m<\kappa( t_{1},\dots,t_{j})$}.
\end{enumerate}
\end{definition}

It is easy to construct a sound function for a finite set $\mathcal{T}$.
For example, let $k=3$ and $\mathcal{T}=\{\vec{t},\vec{s}\}$ where
\[\vec{t}=(t_{1}(x_1,x_3), t_{2}(x_1), t_{3}(x_6))\text{ and }
\vec{s}=(s_{1}(x_1,x_3), s_{2}(x_2,x_9), s_{3}(x_2)).\]
Assuming the terms are all different, a possible sound function for $\mathcal{T}$ is defined as
\begin{align*}
\kappa(t_{1})&=4, \ \ \ \kappa(t_{1}, t_{2})=5, \ \ \ \kappa(t_{1}, t_{2}, t_{3})=7\\
\kappa(s_{1})&=6, \ \ \ \kappa(s_{1}, s_{2})=10, \ \ \kappa(s_{1}, s_{2}, s_{3})=11.
\end{align*}

It is not evident how to construct a sound function for an infinite set $\mathcal{T}$, and clearly it is not possible if $\mathcal{T}$ is uncountable. A very important part of our proof depends on the existence of sound functions for the infinite set $T^k(X)$ when the language is finite. The following lemma guarantees that such a function exists in fact for any countable language.

\begin{lemma}\label{soundlemma} Let $\Len$ be a countable language and $\mathcal{T}\subseteq T^k(X)$ for some $k>0$. Then, there is a sound function for $\mathcal{T}$.
\end{lemma}
\begin{proof}
Let $\beta\colon T(X)\to \N$ be a G\"{o}del numbering for terms in $\Len$ satisfying
\[\text{if } x_m \text{ occurs in }t \text{, then }m<\beta(t)\]
for every term $t$ and every $m\in \N$. This function is injective.

Next, define the function $\kappa\colon \bigcup_{j=1}^kT^{j}(X)\to\N$ by
\[\kappa(t_1,\dots,t_j)=2^{\beta(t_1)}\dots \mathrm{p}_j^{\beta(t_j)}\]
where $\mathrm{p}_j$ is the $j$-th prime number.
It is easy to verify that $\kappa$ is a sound function for $T^{k}(X)$ and, when properly restricted, for any $\mathcal{T}\subseteq T^{k}(X)$.
\end{proof}

Note that, even though the existence of sound functions over the whole $T^k(X)$ for finite languages is required in the proof of the main theorem, the sound functions involved in the statement below are for finite sets.

The overview we presented in the Introduction was given syntactically, closer to classical formulations of the theorem. Here, we opt for a semantic presentation instead. By the Completeness Theorem, both formulations are equivalent, but the semantic approach reflects the nature of our proof and aligns better with potential applications. The statement of Herbrand's Theorem to be proved is as follows.
\begin{theorem}[Herbrand's Theorem]\label{teo1} Let $\K$ be a universal class, $k>0$, and $\varphi(y_1,z_1,\dots,y_k,z_k)$ an open formula. The following are equivalent:
\begin{enumerate}[label=\textup{(\arabic*)}]
\item \label{teoit1}$\K\models \exists y_1\forall z_1\dots\exists y_k\forall z_k \text{ }\varphi(y_1,z_1,\dots,y_k,z_k)$.
\item\label{teoitkappa} There are $\vec{t}_1,\dots, \vec{t}_N\in T^k(X)$ and a sound function $\kappa$ for $\{\vec{t}_1,\dots, \vec{t}_N\}$ such that
\[\K\models\bigvee_{i=1}^N\varphi(t_{i1},x_{\kappa(t_{i1})},\dots,t_{ik},x_{\kappa(t_{i1},\dots,t_{ik})}).\]
\end{enumerate}
Furthermore, $\vec{t}_1,\dots, \vec{t}_N$ and the sound function $\kappa$ can be taken in such a manner that $\{\vec{t}_1,\dots, \vec{t}_N\}\subseteq T^k(x_1,\dots,x_{kN})$, and $\mathrm{Im}(\kappa)\subseteq \{2,\dots, kN+1\}$.
\end{theorem}
The primary significance of the theorem lies in the direction \ref{teoit1} $\Rightarrow$ \ref{teoitkappa}, as it establishes the reduction of the scope of existential quantifiers to a finite set of terms. We will give a model-theoretic proof of this implication.

Regarding the direction \ref{teoitkappa} $\Rightarrow$ \ref{teoit1}, it is quite evident if one understands how sound functions operate, but the proof is notationally and technically complex. For this reason, we define a more compact notation and isolate the central result of that proof in a separate lemma.

Given $k\geq 0$ and a formula $\varphi$, let
\[
\langle\varphi\rangle^{k}_{m}:= \left\{\begin{array}{ll}\exists y_{m+1}\forall z_{m+1}\dots\exists y_k\forall z_k\;
\varphi &\text{ if }0\leq m<k\\
\varphi & \text{ if }m=k.\end{array}\right.\]

Notice that if the free variables of $\varphi$ are in $\{y_1,z_1,\dots,y_k,z_k\}$, then the free variables of $\langle\varphi\rangle^{k}_m$ are in $\{y_1,z_1,\dots,y_m,z_m\}$ for $m>0$, and $\langle\varphi\rangle^{k}_0$ has no free variables.
We will omit the superscript $k$ from the notation as it is fixed hereafter.

In the next lemma, we will apply these formulas to sequences of terms and variables indexed by a sound function. When $m=0$, the expression $\langle\varphi\rangle_{m}(t_{1},x_{\kappa(t_{1})},\dots,t_{m},x_{\kappa(t_{1},\dots,t_{m})})$ is technically incorrect as $\kappa$ is not defined for a  $0$-tuple; by a slight abuse of notation, however, it will be understood simply as the sentence $\langle\varphi\rangle_{0}$.

\begin{lemma}\label{Lemma1}
Let $\varphi(y_1,z_1,\dots,y_k,z_k)$ be a first-order formula, $ \vec{t}_{1},\dots,\vec{t}_{N}\in T^k(X)$, and $\kappa$ a sound function for $\{\vec{t}_{1},\dots,\vec{t}_{N}\}$.
Then, for any $m_1,\dots,m_N\in\{0,\dots,k\}$,
\[\forall x_{1}\dots x_{r} \bigvee_{i=1}^{N} \langle\varphi\rangle_{m_i}(t_{i1},x_{\kappa(t_{i1})},\dots,t_{im_i},x_{\kappa(t_{i1},\dots,t_{im_i})})\vdash \langle\varphi\rangle_{0},\]
where $r=\max\{\kappa(t_{i1},\dots,t_{im_i})\colon m_i>0\}$.
\end{lemma}
\begin{proof}
We assume that no $x_i$ occurs in $\varphi$; this is obvious if $\varphi$ is open. Otherwise, we replace the bound variables $x_i$ with suitable variables.

We will prove this by induction on $M=\sum_{i=1}^N m_i$.
If $M=0$, then $m_i =0$ for each $i\in\{1,\dots,N\}$, and the result trivially follows.
Now suppose $M>0$ and define
\[
I=\{i\colon \ \kappa(t_{i1},\dots,t_{im_i})=r\}
\]
and fix $e\in I$. Note that if $i\in I$, then $\kappa(t_{i1},\dots,t_{im_i})=\kappa(t_{e1},\dots,t_{em_e})=r$ and therefore, by injectivity of sound functions, $m_i =m_e$ and $(t_{i1},\dots,t_{im_{i}})=(t_{e1},\dots,t_{em_{e}})$.

Thus, from $\forall x_{1}\dots x_{r}\bigvee_{i=1}^{N}\langle\varphi\rangle_{m_i }(t_{i1},x_{\kappa(t_{i1})},\dots,t_{im_i},x_{\kappa(t_{i1},\dots,t_{im_i})})$, we obtain
\begin{multline}\label{eq4}
\forall x_{1}\dots x_{r} \ \Bigl(\langle\varphi\rangle_{m_e }(t_{e1},x_{\kappa(t_{e1})},\dots,t_{em_e},x_{r})\lor \\
\bigvee_{i\notin I}\langle\varphi\rangle_{m_i }(t_{i1},x_{\kappa(t_{i1})},\dots,t_{im_i},x_{\kappa(t_{i1},\dots,t_{im_i})})\Bigr).
\end{multline}
Note that, if $i\notin I$, then $x_r$ does not occur in $\langle\varphi\rangle_{m_i}(t_{i1},x_{\kappa(t_{i1})},\dots,t_{im_i},x_{\kappa(t_{i1},\dots,t_{im_i})})$. If $m_i=0$ this is trivial. If $m_i>0$, then $\kappa(t_{i1},\dots,t_{im_i})<r$ and it follows from
\ref{iti} and \ref{itii} of the definition of a sound function. Therefore, (\ref{eq4}) yields
\begin{multline}\label{eq3}
\forall x_{1}\dots x_{r{-}1} \ \Bigl(\bigl(\forall x_r \ \langle\varphi\rangle_{m_e }(t_{e1},x_{\kappa(t_{e1})},\dots,t_{em_e},x_{r})\bigr)\lor \\
\bigvee_{i\notin I}\langle\varphi\rangle_{m_i }(t_{i1},x_{\kappa(t_{i1})},\dots,t_{im_i},x_{\kappa(t_{i1},\dots,t_{im_i})})\Bigr).
\end{multline}
By \ref{iti}, $\kappa(t_{e1},\dots,t_{ej})<r$ for all $j\in \{1,\dots, m_e{-}1\}$ and, by \ref{itii}, $x_{r}$ does not occur in $t_{ej}$ for any $j\in \{1,\dots, m_e\}$. Therefore, $\forall x_{r}\langle\varphi\rangle_{m_e }(t_{e1},x_{\kappa(t_{e1})},\dots,t_{em_e},x_{r})$ implies
\begin{equation}\label{eqinter}
\forall z_{m_e}\langle\varphi\rangle_{m_e }(t_{e1},x_{\kappa(t_{e1})},\dots,t_{em_e},z_{m_e}).
\end{equation}
 Since $t_{em_e}$ is free for $y_{m_e}$ in $\forall z_{m_e }\langle\varphi\rangle_{m_e }(t_{e1},x_{\kappa(t_{e1})},\dots,y_{m_e},z_{m_e})$, (\ref{eqinter}) yields
\[\exists y_{m_e }\forall z_{m_e }\langle\varphi\rangle_{m_e }(t_{e1},x_{\kappa(t_{e1})},\dots,y_{m_e},z_{m_e}).\]
Hence, as $\exists y_{m_e }\forall z_{m_e}\langle\varphi\rangle_{m_e}=\langle\varphi\rangle_{m_e{-}1}$, from (\ref{eq3}), we obtain
\begin{multline}\label{eq5}
\forall x_{1}\dots x_{r{-}1}\Bigl(\langle\varphi\rangle_{m_e{-}1}(t_{e1},x_{\kappa(t_{e1})},\dots,t_{e(m_e{-}1)},x_{\kappa(t_{e1},\dots,t_{e(m_e{-}1)})})\lor \\
\bigvee_{i\notin I}\langle\varphi\rangle_{m_i }(t_{i1},x_{\kappa(t_{i1})},\dots,t_{im_i},x_{\kappa(t_{i1},\dots,t_{im_i})})\Bigr).
\end{multline}

Since $(m_e{-}1)+\sum_{i\notin I}m_i <M$, by the inductive hypothesis we have that (\ref{eq5}) implies $\langle\varphi\rangle_0$, and we have proved Lemma \ref{Lemma1}.
\end{proof}
We are now in a position to prove the main theorem.
\begin{proof}[\textbf{Proof of Theorem \ref{teo1}}]
\
\\
\ref{teoit1} $\Rightarrow$ \ref{teoitkappa}.
Let $\Len$ be the language of $\K$, and let $\Sigma$ be a set of universal axioms for $\K$. Therefore, $\K=\Mod_{\Len}(\Sigma)$ and
\[\Sigma\models \exists y_1\forall z_1\dots\exists y_k\forall z_k \text{ }\varphi(y_1,z_1,\dots,y_k,z_k).\]
By Compactness, there is a finite subset $\Sigma_0\subseteq \Sigma$ such that
\[\Sigma_0\models \exists y_1\forall z_1\dots\exists y_k\forall z_k \text{ }\varphi(y_1,z_1,\dots,y_k,z_k).\]
The language $\Len_0$ containing every symbol in $\Sigma_0\cup\{\varphi\}$ is then clearly finite, and if we prove that $\Sigma_0\models \bigvee_{i=1}^N\varphi(t_{i1},x_{\kappa(t_{i1})},\dots,t_{ik},x_{\kappa(t_{i1},\dots,t_{ik})})$, then the same result holds for $\Sigma$. Therefore, without loss of generality, we assume that $\Len$ is finite.

Let $\kappa$ be a sound function for $T^k(X)$, which exists by Lemma \ref{soundlemma}.

We reserve the notation $T(X)$ for $T_\Len(X)$. Let $\Cons:=\{c_i\colon i\in\N\}$ be a set of distinct constant symbols not in $\Len$, and let $\Len^*:=\Len\cup\Cons$.
Given a term $t(x_{i_1},\dots,x_{i_m})\in T(X)$, we denote $t^*:= t(c_{i_1},\dots, c_{i_m})$. Clearly, every ground term $s\in T_{\Len^*}(\emptyset)$ is a $t^*$ for some $t\in T(X)$.
Now consider an arbitrary $\Len^*$-structure $\A$ such that $\A\models\Sigma$, and let $\B$ be the substructure of $\A$ generated by $\emptyset$, which exists since $\Cons\subseteq\Len^*$, and the elements in $B$ are precisely the interpretations of ground terms of $\Len^*$. As $\Sigma$ is preserved under substructures, we have that $\B\models \Sigma$, and consequently,
\[\B\models\exists y_1\forall z_1\dots\exists y_k\forall z_k \text{ }\varphi(y_1,z_1,\dots,y_k,z_k).\]
It is easy to prove by induction on $m$ that for any $\Len^*$-formula $\psi$ (not necessarily open), if $\B\models\exists y_1\forall z_1\dots\exists y_m\forall z_m\;\psi$, then there exist terms $t_1,\dots, t_m\in T(X)$ such that
\[\B\models \psi(t^*_1,c_{\kappa(t_1)},t^*_2,c_{\kappa(t_1,t_2)},\dots,t^*_m,c_{\kappa(t_1,\dots,t_m)}).\]
In particular, applying this to $m=k$ and $\psi=\varphi$, we obtain
\[\B\models \varphi(t^*_1,c_{\kappa(t_1)},t^*_2,c_{\kappa(t_1,t_2)},\dots,t^*_k,c_{\kappa(t_1,\dots,t_k)}).\]
Since the above is a ground sentence, its satisfaction is preserved under extensions and consequently,
\[\A\models \varphi(t^*_1,c_{\kappa(t_1)},t^*_2,c_{\kappa(t_1,t_2)},\dots,t^*_k,c_{\kappa(t_1,\dots,t_k)}).\]
Because $\A$ is arbitrary, the set of sentences
\[\Sigma\cup\{\neg \varphi(t^*_{1},c_{\kappa(t_{1})},\dots,t^*_{k},c_{\kappa(t_{1},\dots,t_{k})}) \colon \vec{t} \in T^k(X)\}\]
is unsatisfiable. By Compactness, there exist $N>0$ and $\vec{t}_1,\dots,\vec{t}_N\in T^k(X)$ such that
\[\Sigma\models\bigvee_{i=1}^N\varphi(t^*_{i1},c_{\kappa(t_{i1})},t^*_{i2},c_{\kappa(t_{i1},t_{i2})},\dots,t^*_{ik},c_{\kappa(t_{i1},\dots,t_{ik})}).\]
As the constants $c_i$ are not in $\Len$, we conclude that
\begin{equation*}
\Sigma\models\bigvee_{i=1}^N\varphi(t_{i1},x_{\kappa(t_{i1})},t_{i2},x_{\kappa(t_{i1},t_{i2})},\dots,t_{ik},x_{\kappa(t_{i1},\dots,t_{ik})}).
\end{equation*}
Since $\K=\Mod_{\Len}(\Sigma)$, this immediately yields
\[\K\models\bigvee_{i=1}^N\varphi(t_{i1},x_{\kappa(t_{i1})},t_{i2},x_{\kappa(t_{i1},t_{i2})},\dots,t_{ik},x_{\kappa(t_{i1},\dots,t_{ik})}),\]
and we have proved \ref{teoitkappa}.

\ref{teoitkappa} $\Rightarrow$ \ref{teoit1}.
Let $\vec{t}_1,\dots, \vec{t}_N\in T^k(X)$ and let $\kappa$ be a sound function for $\{\vec{t}_1,\dots, \vec{t}_N\}$ such that
\[\K\models\bigvee_{i=1}^N\varphi(t_{i1},x_{\kappa(t_{i1})},t_{i2},x_{\kappa(t_{i1},t_{i2})},\dots,t_{ik},x_{\kappa(t_{i1},\dots,t_{ik})}).\]
The result easily follows from Lemma \ref{Lemma1} by taking $m_i=k$ for every $i$.

To prove the remaining assertion, we perform suitable substitutions.
Without loss of generality, we assume $1\notin \mathrm{Im}(\kappa)$; if this is not the case, we just perform the substitution $x_i\mapsto x_{i+1}$.

Let $I=\mathrm{Im}(\kappa)=\{\kappa(t_{i1},\dots,t_{ij}): 1\leq i\leq N\text{ and } 1\leq j\leq k\}$
and let $M=\max(I)$. Note that all the variables involved in the formula
\begin{equation}\label{psi}\bigvee_{i=1}^N\varphi(t_{i1},x_{\kappa(t_{i1})},t_{i2},x_{\kappa(t_{i1},t_{i2})},\dots,t_{ik},x_{\kappa(t_{i1},\dots,t_{ik})})\end{equation}
belong to the set $\{x_1,\dots, x_M\}$.

The first step is to replace every variable $x_i$ such that $i \notin I$ by $x_1$.
Given a term $t \in T(X)$, we define $t'$ as the result of replacing in $t$ every occurrence of a variable $x_i$ such that $i \notin I$ by $x_1$.
Observe that different terms $t$ and $s$ may result in the same term, i.e., $t' = s'$. Note that the indices of the variables involved in $\{t'_{ij}: 1\leq i\leq N\text{ and } 1\leq j\leq k\}$ are included in $\mathrm{Im}(\kappa)\cup\{1\}$.
With this in mind, let $\lambda$ be defined by
\[\lambda(t'_{i1}, \dots, t'_{ij})=\min\{\kappa(t_{u1}, \dots, t_{uj}): (t'_{u1}, \dots, t'_{uj})=(t'_{i1}, \dots, t'_{ij})\}.\]
It is easy to see that $\lambda$ is a sound function for $\{\vec{t}'_1,\dots,\vec{t}'_N\}$.
Since $\K$ satisfies (\ref{psi}), by instantiation, we have
\begin{equation}\label{line2}
\K \models \bigvee_{i=1}^N\varphi(t'_{i1},x_{\lambda(t'_{i1})},\dots,t'_{ik},x_{\lambda(t'_{i1},\dots,t'_{ik})}).
\end{equation}
Since $\mathrm{Im}(\lambda)\subseteq\mathrm{Im}(\kappa)$, every index of a variable involved in (\ref{line2}) is in $J:=\mathrm{Im}(\kappa)\cup\{1\}$. Now let
 $\alpha\colon\{1,\dots,|J|\}\to J$ be the strictly increasing enumeration of $J$.

Since $|J|\leq kN+1$, by substituting $x_{i}$ with $x_{\alpha^{-1}(i)}$ for every $i\in J$, we obtain the desired result.\end{proof}

\section{Some Applications}
While traditional formulations of Herbrand's theorem rely on a fixed number of quantifier alternations to recursively handle witnesses, our general statement holds uniformly for arbitrary prefixes. To illustrate the structural power of this formulation, we apply it in two corollaries.

Let $\Len$ be a first-order language. We say that a class $\mathcal{Q}$ of $\Len$-structures is a quasivariety if it is a universal class closed under direct products.

The following result provides a characterization for positive formulas satisfied in a quasivariety. We shall consider formulas in prenex disjunctive normal form with an alternating quantifier structure as in the previous section.
\begin{corollary}\label{positive}
Let $\mathcal{Q}$ be a quasivariety and let $\alpha_1,\dots,\alpha_n$ be conjunctions of atomic formulas with variables in $\{y_1,z_1,\dots,y_k,z_k\}$. The following are equivalent:
\begin{enumerate}[label=\textup{(\arabic*)}]
\item\label{positive1} $\mathcal{Q}\models \exists y_{1}\forall z_{1}\dots\exists y_k\forall z_k\; \bigvee_{j=1}^n\alpha_j(y_1,z_1,\dots,y_k,z_k)$.
\item\label{positive2} There are terms $t_1(x_1), t_2(x_1, x_2), \dots, t_k (x_1, \dots, x_k)$ and $j_0\leq n$ such that
\[\mathcal{Q}\models \alpha_{j_0}(t_1,x_2,t_2,x_3,\dots,t_k,x_{k+1}).\]
\end{enumerate}
\end{corollary}

\begin{proof}
\ref{positive1} $\Rightarrow$ \ref{positive2}. From Theorem \ref{teo1}, there are $\vec{t}_1,\dots, \vec{t}_N\in T^k(x_1,\dots,x_{kN})$ and a sound function $\kappa$ for $\{\vec{t}_1,\dots, \vec{t}_N\}$ such that $\mathrm{Im}(\kappa)\subseteq\{2,\dots,kN+1\}$ and
\begin{equation}\label{posit}\mathcal{Q}\models\bigvee_{i=1}^N\bigvee_{j=1}^n\alpha_j(t_{i1},x_{\kappa(t_{i1})},\dots,t_{ik},x_{\kappa(t_{i1},\dots,t_{ik})}).\end{equation}
Let $n_{il}:= \kappa(t_{i1},\dots,t_{il})$ for notational convenience. Clearly $2\leq n_{il}\leq kN+1$.
By a straightforward adaptation of McKinsey's Theorem (see \cite[Theorem 1 and the observation that follows]{McKinsey}) to structures, since $\mathcal{Q}$ is closed under direct products and each $\alpha_j$ is a conjunction of atomic formulas (which behave exactly as atomic formulas with respect to direct products), the disjunction in (\ref{posit}) can be reduced to a single disjunct. Therefore, there are $i_0$ and $j_0$ such that
\[\mathcal{Q}\models \alpha_{j_0}(t_{i_01},x_{n_{i_01}},\dots,t_{i_0k},x_{n_{i_0k}}).\]
Finally, we substitute the variables not in $\{x_{n_{i_0 1}},\dots,x_{n_{i_0 k}}\}$ with $x_1$ and replace every $x_{n_{i_0l}}$ with $x_{l+1}$. This yields the desired result.

\ref{positive2} $\Rightarrow$ \ref{positive1}. This is trivial.
\end{proof}

It is worth noting that the previous corollary generalizes straightforwardly to positive formulas in prenex disjunctive normal form that need not be in an alternating quantifier form starting with an existential quantifier. For example, if the formula is $\forall z_1\exists w_1\exists w_2\forall z_2 \; \bigvee_{i=1}^n\alpha_i(z_1,w_1,w_2,z_2)$
then there are $i_0\leq n$ and terms $t_1(x_1)$ and $t_2(x_1)$ such that \[\mathcal{Q}\models \; \alpha_{i_0}(x_1,t_1,t_2,x_2).\]
This observation will be useful in the proof of the final corollary.

Let $\K$ be a class of structures. An $n$-ary operation $f$ on $\K$ is a choice of an $n$-ary function $f_{\A}\colon A^n\to A$ for each $\A\in \K$. We say that $f$ is \emph{definable in} $\K$ if there is a first-order formula $\varphi(x_1,\dots,x_n,y)$ such that for every $\A\in \K$ and every $a_1,\dots, a_n,b\in A$
\[b=f_{\A}(a_1,\dots,a_n) \text{ if and only if }\A\models \varphi[a_1,\dots,a_n,b].\]
We say that $f$ is a \emph{term operation} on $\K$ if there exists a term $t(x_1,\dots,x_n)$ such that $f_{\A}(a_1,\dots,a_n)=t^{\A}(a_1,\dots,a_n)$ for every $\A\in \K$ and all $a_1,\dots,a_n\in A$.

\begin{corollary}Let $\mathcal{Q}$ be a quasivariety. If $f$ is an operation definable by a positive first-order formula in $\mathcal{Q}$, then $f$ is a term operation on $\mathcal{Q}$.
\end{corollary}

\begin{proof} We will prove the case for unary operations; the generalization is straightforward.
Let $f$ be a unary operation defined by a positive formula $\varphi(v,w)$ in $\mathcal{Q}$. Without loss of generality, we assume
\[\varphi(v,w)= \exists y_1\forall z_1 \dots \exists y_{k}\forall z_{k} \; \bigvee_{j=1}^n\alpha_j(v,w, y_1, z_1,\dots,y_{k},z_{k}),\]
where $\alpha_1,\dots,\alpha_n$ are conjunctions of atomic formulas.

Since $\varphi$ defines an operation, we have that $\mathcal{Q}\models \forall v\exists w \; \varphi(v,w)$, i.e.,
\[\mathcal{Q}\models \forall v\exists w\exists y_1\forall z_1 \dots \exists y_{k}\forall z_{k} \; \bigvee_{j=1}^n\alpha_j(v,w, y_1, z_1,\dots,y_{k},z_{k}).\]
Hence, from Corollary \ref{positive} and the observation following it, we have that there are terms $t_1(x_1), t_2(x_1),t_3(x_1,x_2),\dots,t_{k+1} (x_1, \dots, x_{k})$ and $j_0\leq n$ such that
\[\mathcal{Q}\models \alpha_{j_0}(x_1,t_1,t_2,x_2, t_3, x_3,\dots,t_{k+1},x_{k+1}).\]
This implies
$\mathcal{Q}\models \exists y_1\forall z_1 \dots \exists y_k \forall z_{k}\bigvee_{j=1}^n\alpha_{j}(x_1,t_1, y_1, z_1,\dots,y_{k},z_{k})$ and therefore
\[\mathcal{Q}\models \varphi(x_1,t_1).\]
By the definition of an operation definable by $\varphi$, we conclude that $f_{\A}(a) = t_1^{\A}(a)$ for all $\A \in \mathcal{Q}$ and $a \in A$.
\end{proof}

\section*{Funding}
This work was supported by the Secretaría de Ciencia y Tecnología de la Universidad Nacional de Córdoba, project 33620230100751CB.

\section*{Author Contributions Statement}
The author confirms being the sole contributor of this work and has approved it for publication.

\section*{Competing Interests}
The author declares that there are no competing interests.

\bibliographystyle{plain}
\bibliography{bibliography}

@article{bava2020,
  author    = {Badano, M. and Vaggione, D.},
  title     = {The fundamental theorem of central element theory},
  journal   = {J. Symbolic Logic},
  volume    = {85},
  number    = {4},
  pages     = {1599--1606},
  year      = {2020}
}

@phdthesis{herb,
  author    = {Herbrand, J.},
  title     = {Recherches sur la th{\'e}orie de la d{\'e}monstration},
  school    = {Institut de Math{\'e}matique de l'Universit{\'e} de Paris},
  address   = {Paris},
  year      = {1930},
  note      = {Travaux de l'Institut de Math{\'e}matique de l'Universit{\'e} de Paris}
}

@article{Most,
  author    = {{\L}o{\'s}, J. and Mostowski, A. and Rasiowa, H.},
  title     = {A proof of {H}erbrand's theorem},
  journal   = {J. Math. Pures Appl.},
  volume    = {35},
  number    = {9},
  pages     = {19--24},
  year      = {1956}
}

@article{Most-Fefe,
  author    = {{\L}o{\'s}, J. and Rasiowa, H. and Mostowski, A.},
  title     = {Addition au travail ``{A} proof of {H}erbrand's theorem''},
  journal   = {J. Math. Pures Appl.},
  volume    = {40},
  number    = {9},
  pages     = {129--134},
  year      = {1961}
}

@incollection{Buss,
  author    = {Buss, S. R.},
  title     = {On {H}erbrand's theorem},
  booktitle = {Logic and Computational Complexity},
  editor    = {Leivant, D.},
  series    = {Lecture Notes in Computer Science},
  volume    = {960},
  pages     = {195--209},
  publisher = {Springer},
  year      = {1995},
  doi       = {10.1007/3-540-60178-3_77}
}

@book{Chang,
  author    = {Chang, C. C. and Keisler, H. J.},
  title     = {Model Theory},
  edition   = {3rd},
  publisher = {North-Holland},
  address   = {Amsterdam},
  year      = {1990}
}

@book{Hodges,
  author    = {Hodges, W.},
  title     = {Model Theory},
  publisher = {Cambridge University Press},
  address   = {Cambridge},
  year      = {1993}
}

@article{McKinsey,
  author    = {McKinsey, J. C. C.},
  title     = {The decision problem for some classes of sentences without quantifiers},
  journal   = {J. Symbolic Logic},
  volume    = {8},
  number    = {3},
  pages     = {61--76},
  year      = {1943}
}

@article{Vaggione1999,
  author    = {Vaggione, D.},
  title     = {Central elements in varieties with the {F}raser--{H}orn property},
  journal   = {Adv. Math.},
  volume    = {148},
  number    = {2},
  pages     = {193--202},
  year      = {1999}
}

\end{document}